\documentclass[11pt]{article}

\title{A note on conjectures of F.~Galvin and R.~Rado}
\author{Fran{\c c}ois G. Dorais}
\date{September 13, 2011}
\usepackage{enumerate}
\usepackage{amsmath}
\usepackage{amssymb}
\usepackage{amsthm}

\newcounter{proc}[section]

\newtheorem{theorem}[proc]{Theorem}
\newtheorem{corollary}[proc]{Corollary}
\newtheorem{definition}[proc]{Definition}
\newtheorem{lemma}[proc]{Lemma}
\newtheorem{question}[proc]{Question}

\newcommand{\strong}[1]{\textbf{#1}}
\newcommand{\restrict}{{\upharpoonright}}
\newcommand{\arrows}{{\rightarrow}}
\newcommand{\narrows}{{\nrightarrow}}

\newcommand{\THEN}{\mathrel{\Rightarrow}}
\newcommand{\IFF}{\mathrel{\Leftrightarrow}}
\newcommand{\seqof}[2]{{\langle}#1\,\colon\,#2{\rangle}}
\newcommand{\setof}[2]{{\lbrace}#1\,\colon\,#2{\rbrace}}

\newcommand{\tl}{\mathrel{\triangleleft}}
\newcommand{\tleq}{\mathrel{\trianglelefteq}}
\newcommand{\tg}{\mathrel{\triangleright}}
\newcommand{\tgeq}{\mathrel{\trianglerighteq}}
\newcommand{\edge}{\mathrel{E}}
\newcommand{\dedge}{\mathrel{\vec{E}}}
\newcommand{\goal}[1]{\mathcal{C}_{#1}}

\begin{document}

\maketitle

\begin{abstract}\noindent
  In 1968, Galvin conjectured that an uncountable poset $P$ is the
  union of countably many chains if and only if this is true for every
  subposet $Q \subseteq P$ with size $\aleph_1$.  In 1981, Rado
  formulated a similar conjecture that an uncountable interval graph $G$ is countably
  chromatic if and only if this is true for every induced subgraph $H
  \subseteq G$ with size $\aleph_1$.  Todor{\v c}evi{\' c} has shown
  that Rado's Conjecture is consistent relative to the existence of a
  supercompact cardinal, while the consistency of Galvin's Conjecture
  remains open. In this paper, we survey and collect a variety of
  results related to these two conjectures.  We also show that the
  extension of Rado's conjecture to the class of all chordal graphs is
  relatively consistent with the existence of a supercompact cardinal.
\end{abstract}

\section{Introduction}

Throughout the following, $G$ will denote a (simple loopless) graph
with vertex set $V = V_G$ and edge relation $E = E_G$. For a set $X
\subseteq V$, $G_X$ denotes the induced subgraph with vertex set
$X$. A clique of $G$ is a set $X \subseteq V$ such that $G_X$ is the
complete graph on $X$. Dually, an anticlique is a set $X \subseteq V$
such that $G_X$ is the empty graph on $X$. The conjectures of Galvin
and Rado concern equalities between certain cardinal characteristics
in certain classes of graphs. These cardinal characteristics are the
following.

\begin{itemize}
\item The \strong{clique number} is $$\omega (G) = \sup
  \setof{|X|}{\text{$X$ is a clique of $G$}}.$$
\item The \strong{stability number} is $$\alpha (G) = \sup
  \setof{|X|}{\text{$X$ is an anticlique of $G$}}.$$
\item The \strong{chromatic number} $\chi(G)$ is the smallest size of a cover of $V$ by anticliques.
\item The \strong{clique-cover number} $\theta(G)$ is the smallest
  size of a cover of $V$ by cliques.
\end{itemize}

{\noindent}Clearly, $\omega(G) \leq \chi(G)$ and $\alpha(G) \leq
\theta(G)$. In view of this, it is natural to ask when the equalities
$\omega(G) = \chi(G)$ and $\alpha(G) = \theta(G)$ hold.  It is easy to
check that both equalities fail for the odd cycle $C_{2n+1}$ when $n
\geq 2$.  In 1960, Berge conjectured that the minimal finite graphs for
which these equalities fail are precisely the odd cycles $C_{2n+1}$
and their complements $\overline{C}_{2n+1}$, for $n \geq 2$.  This
fact, the Strong Perfect Graph Theorem, was established by
Chudovsky, Robertson, Seymour, and Thomas in 2002.

Thus, if $G$ is a finite graph that contains no induced copies of the odd cycle $C_{2n+1}$ of
length $2n+1$ nor its complement $\overline{C}_{2n+1}$ for $n \geq 2$,
then the equalities $\omega(G) = \chi(G)$ and $\alpha(G) = \theta(G)$
hold not only for $G$, but also every induced subgraph of $G$. In
fact, we see that the first equality holds for every induced subgraph
of $G$ if and only if the second equality holds for every induced
subgraph of $G$.  This celebrated equivalence, the Perfect Graph
Theorem, was also conjectured by Berge in 1960, and proved by
Lov{\'a}sz in 1972.

\begin{theorem}[Lov{\'a}sz~\cite{Lovasz72b,Lovasz72a};
  Chudnovsky--Robertson--Seymour--Thomas~\cite{ChudnovskyRobertsonSeymourThomas06}]
  \label{thm:perfect}
  The following are equivalent for every graph $G$.
  \begin{enumerate}[\upshape\bfseries(a)]
  \item $G$ contains no induced copies of the odd cycle $C_{2n+1}$ nor
    its complement $\overline{C}_{2n+1}$ for $n \geq 2$.
  \item $\omega(G_X) = \chi(G_X)$ for every finite $X \subseteq V$.
  \item $\alpha(G_X) = \theta(G_X)$ for every finite $X \subseteq V$.
  \item $\alpha(G_X)\omega(G_X) \geq |X|$ for every finite $X
    \subseteq V$.
  \end{enumerate}
  Where $G_X$ denotes the induced subgraph of $G$ with vertex set $X$.
\end{theorem} 

{\noindent}A graph $G$ that satisfies all of these
equivalent properties is known as a \strong{perfect graph}.

Several common types of graphs are known to be perfect. The first to
be identified as such is probably the class of comparability graphs.
Recall that a graph is a \strong{comparability graph} if it has a
transitive orientation or, equivalently, if it is the graph induced by
the comparability relation of a partial ordering of the vertices.

\begin{theorem}[Dilworth~\cite{Dilworth50}]
  Comparability graphs are perfect.
\end{theorem}

{\noindent}Another important class of perfect graphs is the class of
chordal graphs. Recall that a \strong{chordal graph} (also known as a
\strong{triangulated graph}) is a graph that has no induced copies of
the cycle $C_n$ for $n \geq 4$. (See
Theorem~\ref{thm:perfectelimination} for an alternate
characterization.)

\begin{theorem}[Hajnal--Sur{\'a}nyi~\cite{HajnalSuranyi58}; Berge~\cite{Berge60}]
  Chordal graphs are perfect.
\end{theorem}

{\noindent}\strong{Interval graphs} (i.e., intersection graphs of
families of non-empty convex subsets of a linear order) are also
perfect. This can be seen in two ways: because every interval graph is
a chordal graph, or because the complement of every interval graph is
a comparability graph.

The following result is a typical use of the compactness theorem in
graph theory.

\begin{theorem}
  Let $G = (V,E)$ be a graph and let $k$ be a positive integer.
  \begin{enumerate}[\upshape\bfseries(a)]
  \item $\chi(G) \leq k$ if and only if $\chi(G_X) \leq k$ for every
    finite $X \subseteq V$.
  \item $\theta(G) \leq k$ if and only if $\theta(G_X) \leq k$ for
    every finite $X \subseteq V$.
  \end{enumerate}
\end{theorem}

{\noindent}For perfect graphs, we have a very strong form of this
fact.

\begin{corollary}\label{cor:finperfect}
  Let $G$ be a perfect graph and let $k$ be a positive integer.
  \begin{enumerate}[\upshape\bfseries(a)]
  \item\label{cor:finperfect:col} $\chi (G) \leq k$ if and only if
    $\chi (G_X) \leq k$ for every $X \in [V]^{k + 1}$.
  \item\label{cor:finperfect:cov} $\theta (G) \leq k$ if and only if
    $\theta (G_X) \leq k$ for every $X \in [V]^{k + 1}$.
  \end{enumerate}
\end{corollary}

\begin{proof}
  It is enough to prove~\eqref{cor:finperfect:col}
  since~\eqref{cor:finperfect:cov} is dual. Note that $\chi (G_X) \leq
  k$ for every $X \in [V]^{k + 1}$ if and only if $\omega (G) \leq
  k$. By Theorem~\ref{thm:perfect}, $\omega (G) \leq k$ if and only if
  $\chi (G_X) \leq k$ for every finite $X \subseteq V$.
\end{proof}

{\noindent}It is natural to ask whether the same holds if one replaces
$k$ by an infinite cardinal $\kappa$ and $k+1$ by it cardinal
successor $\kappa^+$. We will concentrate in the first case, $\kappa =
\aleph_0$ and $\kappa^+ = \aleph_1$.

\begin{definition}
  Let $\Gamma$ be a class of graphs. We use $\goal{\chi}$ and
  $\goal{\theta}$ denote the following dual statements.
  \begin{description}
  \item[$(\goal{\chi})$] For every $G \in \Gamma$, $\chi (G) \leq
    \aleph_0$ if and only if $\chi (G_X) \leq \aleph_0$ for every $X
    \in [V]^{\aleph_1}$.
  \item[$(\goal{\theta})$] For every $G \in \Gamma$, $\theta (G) \leq
    \aleph_0$ if and only if $\theta (G_X) \leq \aleph_0$ for every $X
    \in [V]^{\aleph_1}$.
  \end{description}
\end{definition}

{\noindent}In 1968, Galvin~\cite{Galvin08} conjectured that $\goal{\theta}$ holds for
the class of comparability graphs.%
In 1981, Rado~\cite{Rado81} conjectured that the
class of interval graphs has property $\goal{\chi}$. The consistency,
relative to the existence of a supercompact cardinal, of Rado's
Conjecture was then established by Todor{\v
  c}evi{\'c}~\cite{Todorcevic83} in 1983. In~\cite{Todorcevic83}
and~\cite{Todorcevic91}, Todor{\v c}evi{\' c} shows that large
cardinals are indeed necessary to establish the consistency of Rado's
Conjecture.

In this paper, we will show that Todor{\v c}evi{\'c}'s result on the
consistency of Rado's Conjecture can be extended to the consistency of
$\goal{\chi}$ for the class of all chordal graphs.

\begin{theorem}\label{thm:superrado}
  Each of the following statements implies the next.
  \begin{enumerate}[\upshape\bfseries(a)]
  \item \label{thm:superrado:treeable} $\goal{\chi}$ holds for the
    class of $\sigma$-treeable graphs.
  \item \label{thm:superrado:chordal} $\goal{\chi}$ holds for the
    class of chordal graphs.
  \item \label{thm:superrado:interval} $\goal{\chi}$ holds for the
    class of interval graphs \textup(Rado's Conjecture\textup).
  \end{enumerate}
  Furthermore, these statements are all consistent relative to the
  existence of a supercompact cardinal.
\end{theorem}

{\noindent}This theorem will be proved in Section~\ref{sec:superrado}
(where we also define $\sigma$-treeable graphs). We do not know if any
of the implications of Theorem~\ref{thm:superrado} are strict since
the same technique is used to prove the consistency in all cases.

We will also provide a proof of the following result of Todor{\v
  c}evi{\'c} which shows that Rado's Conjecture is equivalent to the
restriction of Galvin's Conjecture to the class of finite-dimensional
comparability graphs.

\begin{theorem}[Todor{\v c}evi{\'c}]\label{thm:radoequiv}
  The following are equivalent.
  \begin{enumerate}[\upshape\bfseries(a)]
  \item \label{thm:radoequiv:interval} $\goal{\chi}$ holds for the
    class of interval graphs \textup(Rado's Conjecture\textup).
  \item \label{thm:radoequiv:twodimanti} $\goal{\chi}$ holds for the
    class of $2$-dimensional comparability graphs.
  \item \label{thm:radoequiv:twodimchain} $\goal{\theta}$ holds for
    the class of $2$-dimensional comparability graphs.
  \item \label{thm:radoequiv:findimchain} $\goal{\theta}$ holds for
    the class of finite-dimensional comparability graphs.
  \end{enumerate}
\end{theorem}

{\noindent}The equivalence of~\eqref{thm:radoequiv:interval}
and~\eqref{thm:radoequiv:findimchain} appears without proof
in~\cite[Remark~4.6]{Todorcevic11}.  A proof of this theorem will be
provided in Section~\ref{sec:radoequiv} (where we also define
$n$-dimensional comparability graphs).

While the consistency of Galvin's Conjecture remains open, the above
results lead us to the following more general question.

\begin{question}
  Is it consistent, relative to large cardinals, that $\goal{\chi}$
  and, equivalently, $\goal{\theta}$ hold for the class of perfect
  graphs?
\end{question}

In view of Theorem~\ref{thm:superrado}, it is natural to ask about
$\goal{\theta}$ for the class of chordal graphs.  It turns out that
$\goal{\theta}$ is simply true for this class.  In fact, property
$\goal{\theta}$ holds for the broader class of \strong{squarefree}
graphs, i.e., graph that do not contain induced copies of the square
$C_4$. This follows from a result of Wagon.

\begin{theorem}[Wagon~\cite{Wagon78}]
  Suppose $G$ is a squarefree graph such that $\alpha(G) \leq
  \aleph_0$.  Then $\theta(G) > \aleph_0$ if and only if $G$ contains
  an induced copy of the comparability graph of a Suslin tree.
\end{theorem}

{\noindent}Since Suslin trees have size $\aleph_1$, we have the
following immediate corollary.

\begin{corollary}\label{cor:wagon}
  $\goal{\theta}$ holds for the class of squarefree graphs, and hence
  for the class of chordal graphs.
\end{corollary}

{\noindent}The techniques used by Wagon suggest that many squarefree
graphs are $\sigma$-treeable, so there is a chance that the dual of
Corollary~\ref{cor:wagon} is consistent relative to large cardinals.

\begin{question}
  Is it consistent, relative to large cardinals, that $\goal{\chi}$
  holds for the class of squarefree graphs?
\end{question}

\section{Results of Abraham and Todor{\v c}evi{\' c}}

In this section we summarize some earlier theorems that shed some
light on the conjectures of Galvin and Rado.  The first due to
Abraham and the second due to Todor{\v c}evi{\' c}.  Abraham's
result shows that Galvin's Conjecture holds for the class of
comparability graphs without infinite anticliques. Todor{\v c}evi{\'
  c}'s result gives several equivalent forms of Rado's Conjecture in
terms of one-dimensional partition relations for posets.

In 1963, Perles~\cite{Perles63} showed that Dilworth's Theorem
($\alpha(G) = \theta(G)$ for finite comparability graphs) fails for
infinite comparability graphs by observing that the cartesian product
$\omega_1 \times \omega_1$ has no infinite antichains but cannot be
covered by countably many chains.  This example can be generalized as
follows.

\begin{definition}[Abraham~\cite{Abraham87}]
  A poset $P$ is of \strong{Perles type} if there is an enumeration
  $\seqof{p_{\alpha}}{\alpha < \omega_1}$ of $P$ and a function $f :
  \omega_1 \to \omega_1$ such that $|f^{-1}(\alpha)| = \aleph_1$ for
  every $\alpha < \omega_1$, and $\alpha < \beta \land f(\alpha) >
  f(\beta)$ imply that $p_{\alpha}$ and $p_{\beta}$ are incomparable.
\end{definition}

{\noindent}This definition captures the essential features of
$\omega_1 \times \omega_1$ that were used in Perles's counterexample.
Abraham~\cite{Abraham87} then showed that these are essentially the
only counterexamples to Dilworth's Theorem that don't have infinite
antichains.

\begin{theorem}[Abraham~\cite{Abraham87}]\label{thm:abraham}
  Suppose $P$ is a poset without infinite antichains. Then $P$ is the
  union of countably many chains if and only if it does not contain a
  poset of Perles type.
\end{theorem}

{\noindent}Since the posets of Perles type all have size $\aleph_1$,
it follows immediately that:

\begin{corollary}\label{cor:abraham}
  $\goal{\theta}$ holds for the class comparability graphs without
  infinite anticliques.
\end{corollary}

To state Todor{\v c}evi{\' c}'s result, it is convenient to introduce
some ``Hungarian notation'' for one-dimensional partitions of posets.
If $\psi$ is a partial order type and $\kappa$ is a cardinal, we
write $P \arrows (\psi)^1_\kappa$ if for every coloring
$c:P\to\kappa$, there is a $Q \subseteq P$ with order type $\psi$ such
that $c$ is constant on $Q$; $P \narrows (\psi)^1_\kappa$ denotes the
negation of this statement. Generalizing this notation a little, if
$\psi_1,\dots,\psi_k$ are partial order types and $\kappa$ is a
cardinal, we write $P \arrows (\psi_1 \vee \cdots \vee
\psi_k)^1_\kappa$ if for every coloring $c:P\to\kappa$, there is a $Q
\subseteq P$, with order type among $\psi_1,\dots,\psi_k$, such that
$c$ is constant on $Q$; again $P \narrows
(\psi_1\vee\cdots\vee\psi_k)^1_\kappa$ denotes the negation of this
statement.

We will mostly be interested in the negative cases when $\psi \in
\{2,\omega,\omega^*\}$.  Indeed, $P \narrows (2)^1_\kappa$ simply
means that $P$ is the union of at most $\kappa$ antichains, i.e.,
$\theta(G_P) \leq \kappa$ where $G_P$ is the comparability graph of
$P$.  Similarly, $\phi \narrows (\omega^*)^1_\kappa$ (resp.\
$\phi\narrows(\omega)^1_\kappa$) means that $P$ is the union of at
most $\kappa$ well-founded (resp.\ conversely well-founded) subsets.
Finally, $\phi \narrows (\omega \vee \omega^*)^1_\kappa$ means that
$P$ is the union of $\kappa$ subsets without infinite chains.

\begin{theorem}[Todor{\v c}evi{\' c}~\cite{Todorcevic83}]
  \label{thm:todorcevic}
  The following are equivalent to Rado's Conjecture \textup(i.e.,
  $\goal{\chi}$ holds for the class of interval graphs\textup).
  \begin{enumerate}[\upshape\bfseries(a)]
  \item \label{thm:todorcevic:treespecial}For every tree $T$, $T
    \narrows (2)^1_{\omega}$ if and only if $U \narrows
    (2)^1_{\omega}$ for every $U \in [T]^{\aleph_1}$.
  \item \label{thm:todorcevic:special}For every poset $P$, $P \narrows
    (\omega)^1_{\omega}$ if and only if $Q \narrows
    (\omega)^1_{\omega}$ for every $Q \in [P]^{\aleph_1}$.
  \item \label{thm:todorcevic:dualspecial}For every poset $P$, $P
    \narrows (\omega^{\ast})^1_{\omega}$ if and only if $Q \narrows
    (\omega^{\ast})^1_{\omega}$ for every $Q \in [P]^{\aleph_1}$.
  \item \label{thm:todorcevic:bispecial}For every poset $P$, $P
    \narrows (\omega \lor \omega^{\ast})^1_{\omega}$ if and only if $Q
    \narrows (\omega \lor \omega^{\ast})^1_{\omega}$ for every $Q \in
    [P]^{\aleph_1}$.
  \end{enumerate}
\end{theorem}

\begin{proof}
  The equivalence of Rado's Conjecture with
  (\ref{thm:todorcevic:treespecial}) and
  (\ref{thm:todorcevic:special}) is \cite[Theorem 6]{Todorcevic83};
  (\ref{thm:todorcevic:dualspecial}) is equivalent to
  (\ref{thm:todorcevic:special}), by duality;
  (\ref{thm:todorcevic:bispecial}) follows from the combination of
  (\ref{thm:todorcevic:special}) and
  (\ref{thm:todorcevic:dualspecial}); (\ref{thm:todorcevic:bispecial})
  implies (\ref{thm:todorcevic:treespecial}) since $T \narrows
  (2)^1_{\omega}$, $T \narrows (\omega)^1_{\omega}$, and $T \narrows
  (\omega \lor \omega^{\ast})^1_{\omega}$ are all equivalent for a
  tree $T$. % ???Who???
\end{proof}

\section{Finite-Dimensional Comparability Graphs}\label{sec:radoequiv}

A \strong{$n$-dimensional poset} is a poset $P$ whose order relation
is the intersection of $n$ linear orders, i.e., if there are linear
orders ${\leq_1}, \ldots, {\leq_n}$ on the points of $P$ such that $x
\leq_P y \IFF x \leq_1 y \land \cdots \land x \leq_n y$.  It turns out
that the dimension of a poset is an invariant of its comparability
graph.

\begin{theorem}[Trotter--Moore--Sumner~\cite{TrotterMooreSumner76}]
  If the graph $G$ has an $n$-dimensional transitive orientation, then
  every transitive orientation of $G$ is $n$-dimensional.
\end{theorem}

{\noindent}Thus, it makes sense to say that $G$ is a
\strong{$n$-dimensional comparability graph} if $G$ has an
$n$-dimensional transitive orientation.

The class of $2$-dimensional comparability graphs is especially
interesting since it is self-dual.

\begin{theorem}[Pnueli--Lempel--Even~\cite{PnueliLempelEven71}]
  A graph $G$ is a $2$-dimensional comparability graph if and only if
  $G$ and its complement $\overline{G}$ are both comparability graphs.
\end{theorem}

{\noindent}This last result immediately implies the equivalence of
(\ref{thm:radoequiv:twodimanti}) and (\ref{thm:radoequiv:twodimchain})
in Theorem~\ref{thm:radoequiv}. The next result shows that
(\ref{thm:radoequiv:twodimchain}) implies
(\ref{thm:radoequiv:findimchain}) in Theorem~\ref{thm:radoequiv}. This
establishes the equivalence of the last three statements of
Theorem~\ref{thm:radoequiv}.

\begin{theorem}
  If $\goal{\theta}$ holds for the class of $2$-dimensional
  comparability graphs, then $\goal{\theta}$ holds for the class of
  finite-dimensional comparability graphs.
\end{theorem}

\begin{proof}
  We proceed by induction on dimension. Suppose that $\goal{\theta}$
  holds for every $n$-dimensional comparability graph. Let $P =
  (V,{\leq})$ be a $(n + 1)$-dimensional poset. Then there are a
  $n$-dimensional partial order ${\leq_0}$ and a linear order
  ${\leq_1}$ on $V$ such that $u \leq v \IFF u \leq_0 v \land u \leq_1
  v$. Write $P_0 = (V,{\leq_0})$ and $P_1 = (V,{\leq_1})$. If every $U
  \in [V]^{\aleph_1}$ is the union of countably many ${\leq}$-chains,
  then it is also the union of countably many
  ${\leq_0}$-chains. Therefore, by the induction hypothesis, $V$ is
  the union of countably many ${\leq_0}$-chains, say $V = \bigcup_{n =
    0}^{\infty} C_n$ where each $C_n$ is a ${\leq_0}$-chain. Now the
  restriction of ${\leq}$ to $C_n$ is $2$-dimensional as ${\leq_0}$
  and ${\leq_1}$ are both linear orders on $C_n$. Also, by hypothesis,
  every $D \in [C_n]^{\aleph_1}$ is the union of countably many
  ${\leq}$-chains. Since $\goal{\theta}$ holds for $2$-dimensional
  comparability graphs, each $C_n$ is itself the union of countably
  many ${\leq}$-chains. Gathering these smaller chains together, we
  find that $V$ is the union of countably many ${\leq}$-chains.
\end{proof}

For the last equivalent form of Theorem~\ref{thm:radoequiv}, we appeal
to the results of Todor{\v c}evi{\' c} and Abraham from the previous
section.

\begin{theorem}
  Rado's Conjecture implies that $\goal{\theta}$ holds for the class of
  $2$-dimensional comparability graphs.
\end{theorem}

\begin{proof}
  Let $P$ be a $2$-dimensional poset and let $P'$ be a poset whose
  comparability graph is the complement of that of $P$. Assume that
  every $Q \in [P]^{\aleph_1}$ is the union of countably many chains
  or, dually, every $Q' \in [P']^{\aleph_1}$ is the union of countably
  many antichains. Then we have $Q' \narrows (\omega \lor
  \omega^{\ast})^1_{\omega}$ (indeed $Q' \narrows (2)^1_{\omega}$) for
  every $Q' \in [P']^{\aleph_1}$. Hence, $P' \narrows (\omega \lor
  \omega^{\ast})^1_{\omega}$, by Theorem~\ref{thm:todorcevic}.  Thus,
  $P'$ is the union of countably many sets each of which has no
  infinite chains. It follows by duality that $P = \bigcup_{n =
    0}^{\infty} R_n$ where each $R_n$ has no infinite antichains. Now,
  every $Q \in [R_n]^{\aleph_1} \subseteq [P]^{\aleph_1}$ is the
  union of countably many chains. It follows from
  Corollary~\ref{cor:abraham}, that each $R_n$ is the union of
  countably many chains. Gathering these chains together, we
  see that $P$ is the union of countably many chains.
\end{proof}

{\noindent}This shows that Rado's Conjecture implies Galvin's
Conjecture for $2$-dimensional posets. For the converse, we show that
Galvin's Conjecture for $2$-dimensional posets implies the first
equivalent form of Rado's Conjecture in Theorem~\ref{thm:todorcevic}.

\begin{theorem}
  If $\goal{\chi}$ holds for the class of $2$-dimensional
  comparability graphs then, for every tree $T$, we have $T \narrows
  (2)^1_{\omega}$ if and only if $U \narrows (2)^1_{\omega}$ for every
  $U \in [T]^{\aleph_1}$.
\end{theorem}

\begin{proof}
  It suffices to observe that every tree $T$ is a $2$-dimensional
  poset, which can be seen by lexicographically ordering $T$ in two
  opposite ways.
\end{proof}

\section{Interval, Chordal, and $\sigma$-Treeable
  Graphs}\label{sec:superrado}

The following characterization of chordal graphs is due to
Fulkerson and Gross~\cite{FulkersonGross65} in the finite
case; the infinite case follows by a simple application of the
Compactness Theorem.  An orientation ${\dedge}$ of $G = (V,E)$ is said
to be a \strong{simplicial orientation} if it is acyclic and $S_v =
\setof{u \in V}{u \dedge v}$ is a clique in $G$ for every $v \in V$.

\begin{theorem}[Fulkerson--Gross~\cite{FulkersonGross65}]
  \label{thm:perfectelimination}
  A graph $G = (V, E)$ is chordal if and only if it has a simplicial
  orientation.
\end{theorem}

{\noindent}With this result, it is easy to show that interval graphs
are chordal.

\begin{corollary}
  Every interval graph is chordal.
\end{corollary}

\begin{proof}
  Let $G = (V, E)$ be an interval graph as witnessed by the family of
  intervals $\seqof{I_v}{v \in V}$ of a linear order $L$.  Define $u
  \dedge v$ iff $I_u \cap I_v$ is a nonempty initial subinterval of
  $I_v$.  (If some of the intervals are equal, break ties using a
  linear ordering of $V$.) It is easy to check that ${\dedge}$ is a
  simplicial orientation of $V$.
\end{proof}

{\noindent}It follows immediately that (\ref{thm:superrado:chordal})
implies (\ref{thm:superrado:interval}) in Theorem~\ref{thm:superrado}.

Before we define the class of $\sigma$-treeable graphs, let us make an
observation to motivate the definition.  A poset $R = (V,{\tl})$ is
\strong{ramified} if every initial interval $R[{\tl v}]$ is linearly
ordered by ${\tl}$ for each $v \in V$.  Thus a tree is simply a
well-founded ramified poset.

\begin{theorem}\label{thm:chordalramified}
  If $G = (V,E)$ is a chordal graph, then the transitive closure of
  any simplicial orientation of $G$ is a ramified ordering of $V$.
\end{theorem}

\begin{proof}
  Let ${\tl}$ be the transitive closure of a simplicial orientation
  ${\dedge}$ of $G$. For $v \in V$, let $S_v = \setof{w \in V}{w
    \dedge v}$. Then define $S^0_v = \{v\}$ and $S^{n+1}_v = \bigcup
  \setof{S_w}{w \in S^n_v}$. Note that $S^1_v = S_v$ and $w \tleq v$
  iff $w \in \bigcup_{n=0}^\infty S^n_v$.
  
  We want to show that if $u, v \tleq w$ then $u \tleq v$ or $u \tgeq
  v$. We proceed by induction on $m$ where $u \in S^m_w$.
  
  For $m = 0$, we have $u = w$ and hence $u \tgeq v$.
  
  For $m = 1$, let $v = v_0 \dedge v_1 \dedge \cdots \dedge v_n = w$
  witness that $v \tleq w$. Let $p = \min \setof{i}{u \tl v_i}$.  Note
  that $u \dedge v_i$ for $i = p, \ldots, n$. (This is clear for $i =
  n$ since $u \edge w$ by definition of $S_w$. Suppose that $u \dedge
  v_{i + 1}$ and $i \geq p$, then $u, v_i \in S_{v_{i + 1}}$, which
  means that $u \dedge v_i$ since $S_{v_{i + 1}}$ is a clique and $u
  \tl v_i$.) If $p = 0$ then it follows immediately that $u \tl v_0 =
  v$. If $p > 0$, then note that $v_{p - 1} \dedge u$ or $v_{p-1} = u$
  since $u, v_{p - 1} \in S_{u_p}$, $S_{u_p}$ is a clique, and $u
  \not\dedge v_{p - 1}$.  Therefore, $v = v_0 \tleq v_{p - 1} \tl u$.
  
  For $m > 1$, note that $u \in S_x$ for some $x \in S^{m-1}_w$. By
  the induction hypothesis, either $x \tl v$, $x = v$, or $x \tg
  v$. If $x \tleq v$, then $u \tl v$ by transitivity of ${\tl}$. If $x
  \tg v$, then the result follows from the case $m = 1$, since $u \in
  S_x = S^1_x$.
\end{proof}

A graph $G = (V, E)$ is \strong{$\sigma$-treeable} if it is contained
in the comparability graph of a ramified ordering ${\tl}$ of $V$ which
has the additional property that $|V[{\tl v}]| \leq \aleph_0$ for
every $v \in V$.  The next lemma will perhaps clarify our choice of
terminology.

\begin{lemma}\label{lem:galvin}
  If $G = (V, E)$ is $\sigma$-treeable, as witnessed by the ramified
  ordering ${\tl}$ of $V$, then there is a partition $V = \bigcup_{n =
    0}^{\infty} V_n$ such that the restriction of ${\tl}$ to each
  $V_n$ is a tree of height at most $\omega_1$.
\end{lemma}

\begin{proof} (Due to Galvin, cf.~\cite{Todorcevic83}.) Fix a
  well-ordering ${\prec}$ of $V$. For each $v \in V$, let $f_v :
  V[{\tleq v}] \to \omega$ be an injection. Define, $f : V \to \omega$
  by $f(u) = f_v(u)$ where $v$ is the ${\prec}$-first element of $V$
  such that $u \tleq v$. We claim that the restriction of ${\tl}$ to
  each $V_n = f^{- 1}(n)$ is well-founded.
  
  Suppose that $u_0 \tgeq u_1 \tgeq \cdots$ is a descending sequence
  of elements of $V_n$. Let $v_i$ be the ${\prec}$-first $v \in V$
  with $u_i \tleq v$. Note that $f_{v_i}(u_i) = f(u_i) = n$ for each
  $i < \omega$.  Note also that $v_0 \succeq v_1 \succeq \cdots$ Since
  ${\prec}$ is a well-ordering, there are $v$ and $k$ such that $v_i =
  v$ for $i \geq k$.  Since $f_v$ is an injection we have $u_i =
  f_v^{- 1}(n)$ for $i \geq k$.  Thus $u_0 \tgeq u_1 \tgeq \cdots$ is
  eventually constant, which shows that ${\tl}$ is well-founded on
  $V_n$.
\end{proof}

{\noindent}If $G = (V,E)$ is any graph such that $\chi(G_X) \leq
\aleph_0$ for every $X \in [V]^{\aleph_1}$, then we certainly have
$\omega(G) \leq \aleph_0$.  If, moreover, $G$ is chordal and
${\dedge}$ is a simplicial orientation of $G$, then $|S_v| \leq
\aleph_0$ for each $v \in V$.  It then follows that $|V[{\tl v}]| \leq
\aleph_0$ for each $v \in V$ where ${\tl}$ is the transitive closure
of ${\dedge}$.  Therefore, every chordal graph such that $\chi(G_X)
\leq \aleph_0$ for every $X \in [V]^{\aleph_1}$, is $\sigma$-treeable.
This shows that (\ref{thm:superrado:treeable}) implies
(\ref{thm:superrado:chordal}) in Theorem~\ref{thm:superrado}.

To complete the proof of Theorem~\ref{thm:superrado}, it remains to
prove the consistency of $\goal{\chi}$ for the class of
$\sigma$-treeable graphs, relative to the existence of a supercompact
cardinal.  Rather than giving a forcing proof the consistency of
$\goal{\chi}$ for $\sigma$-treeable graphs, as in~\cite{Todorcevic83},
we will use the Global Game Reflection Principle ($\mathsf{GRP}^+$)
of~\cite{Konig04}. Let $\mathcal{S} \subseteq (A \times B)^{<
  \omega_1}$ be a tree and let $[\mathcal{S}] = \setof{s \in (A \times
  B)^{\omega_1}}{(\forall \alpha < \omega_1) (s \restrict \alpha \in
  \mathcal{S})}$. Consider a two player game $\mathbb{G}(\mathcal{S})$
of length $\omega_1$ where in each round $\alpha < \omega_1$, Player~I
plays $a_{\alpha} \in A$, Player~II responds with $b_{\alpha} \in B$,
and Player~II wins if $\seqof{(a_{\alpha}, b_{\alpha})}{\alpha <
  \omega_1} \in [\mathcal{S}]$. If $X \subseteq A$, then the
restricted game $\mathbb{G}(\mathcal{S}|X)$ is defined similarly
except that Player~I can only play elements of $X$.
\begin{description}
\item[$(\mathsf{GRP}^+)$] If $\mathcal{S} \subseteq (A \times B)^{<
    \omega_1}$ is a tree, $\mathcal{C} \subseteq [A]^{\aleph_1}$ is an
  $\omega_1$-club, and Player~II has a winning strategy in the
  restricted game $\mathbb{G}(\mathcal{S}|X)$ for every $X \in
  \mathcal{C}$, then Player~II has a winning strategy in the
  unrestricted game $\mathbb{G}(\mathcal{S})$.
\end{description}
It is known that this principle has considerable large cardinal
strength, but no more than a supercompact cardinal.  In fact, the
consistency of $\mathsf{GRP}^+$ can be obtained by the L{\'e}vy
collapse of a supercompact cardinal to $\aleph_2$.

\begin{theorem}[K{\"o}nig~\cite{Konig04}] If $\kappa$ is supercompact,
  then $\mathsf{Coll}(\aleph_1,{<}\kappa) \Vdash \mathsf{GRP}^+$.
\end{theorem}

{\noindent}It is also observed in~\cite{Konig04} that Rado's Conjecture
follows from $\mathsf{GRP}^+$. Here we prove the more general result that
$\mathsf{GRP}^+$ implies $\goal{\chi} $ for $\sigma$-treeable graphs.

\begin{theorem}
  $\mathsf{GRP}^+$ implies that $\goal{\chi} $ holds for
  $\sigma$-treeable graphs.
\end{theorem}

\begin{proof}
  Let $G = (V, E)$ be a $\sigma$-treeable graph as witnessed by the
  partial ordering ${\tl}$ of $V$. By Lemma~\ref{lem:galvin}, we may
  assume that $(V,{\tl})$ is a tree of height at most $\omega_1$.
  
  Consider the game $\mathbb{G}_{\chi}(G)$ of length $\omega_1$ where,
  in each round $\alpha < \omega_1$, Player~I plays $v_{\alpha} \in
  V$, Player~II responds with $c_{\alpha} \in \omega$, and Player~II
  wins iff $v_{\alpha} \edge v_{\beta} \THEN c_{\alpha} \neq
  c_{\beta}$ for all $\alpha < \beta < \omega_1$. The fact that
  $\chi(G_W) \leq \aleph_0$ for every $W \in [V]^{\aleph_1}$ clearly
  implies that Player~II has a winning strategy for the restricted
  game $\mathbb{G}_{\chi}(G|W)$. Therefore, by $\mathsf{GRP}^+$,
  Player~II has a winning strategy in the unrestricted game
  $\mathbb{G}_{\chi}(G)$.
  
  Define the coloring $c : V \to \omega$ as follows. Suppose that $v
  \in V$ has height $\eta < \omega_1$ and let
  $\seqof{v_{\alpha}}{\alpha \leq \eta}$ enumerate the branch
  $V[{\tleq v}]$ in ${\tl}$-order (so $v_{\eta} = v$). Consider the
  sequence $\seqof{v_{\alpha}}{\alpha \leq \eta}$ as a sequence of
  moves for Player~I in the game $\mathbb{G}_{\chi}(G)$ and let
  $\seqof{c_{\alpha}}{\alpha \leq \eta}$ be the sequence of Player~II
  responses according to her winning strategy. Then set $c(v) =
  c_{\eta}$.  Note that $c_{\alpha} = c(v_{\alpha})$ for every $\alpha
  \leq \eta$. Since Player~II was using her winning strategy in this
  play, it follows that $v_{\alpha} \edge v_{\eta} \THEN c(v_{\alpha})
  = c_{\alpha} \neq c_{\eta} = c (v_{\eta})$. Therefore, $c : V \to
  \omega$ is a proper coloring of $G$ and hence $\chi(G) \leq
  \aleph_0$.
\end{proof}

\section*{Acknowledgements}
I would like to thank professor Fred Galvin for sharing with me the
history of his conjecture. I would also like to thank professor
Justin Moore for many valuable discussions and for suggesting this line
of inquiry.
I am also indebted to the anonymous referees for their valuable
comments and suggestions.

\bibliographystyle{amsplain}
\bibliography{galvin}

\end{document}